\documentclass[reqno]{amsart}

\usepackage{amssymb,amsmath,amsthm,bbm,enumerate,slashed}
\usepackage[margin=1in]{geometry}

\DeclareMathOperator\Div{div}
\DeclareMathOperator\tr{tr}

\newcommand{\cS}[1]{\mathcal{S}{[#1]}}

\newtheorem{thm}{Theorem}
\newtheorem{prop}[thm]{Proposition}

\newtheorem{lem}[thm]{Lemma}

\newtheorem*{thm*}{Theorem}
\newtheorem*{prop*}{Proposition}

\title{On evolution of mean curvature and minimal areas in general codimension for inverse problems} 
\author{Chun-Kai Kevin Chien}
\email{ck.chien@uni-muenster.de}

\begin{document}

\begin{abstract}
Many seemingly standard formulas for the evolution of the mean curvature vector in higher codimension assume only normal variations, which are ill-suited for analyzing the evolution of submanifolds-with-boundary whose boundaries may not be moving orthogonally at each time. Such situations arise when one considers inverse problems for minimal areas, for example, where such boundaries are confined to a larger ambient boundary.

In this note, we provide coordinate-free identities for use in such problems, including higher order derivatives of mean curvature. We translate these identities to some of the existing inverse problem literature in codimension one. We also hope readers may find the computational framework helpful in other situations arising in geometric analysis.
\end{abstract}

\maketitle

\section*{Statement on AI Use}
No AI was used in the writing, literature review, ideas, or calculations contained in this document except for the singular use of Claude 4.5 to verify the statement of Lemma \ref{lem:divnormtotan} in the special case of $Y=\partial_tF$.

\section{Introduction}
Let $F:I^k\times \Sigma^m\to (M^n,g)$ be a smooth $k$-parameter family of immersions of a smooth, orientable, compact manifold-with-boundary $\Sigma$ into a (pseudo-)Riemannian manifold-with-boundary $(M,g)$, with parameter $t=(t^1,\ldots, t^k)$-dependent induced (Riemannian) metric $\sigma(t):=F_t^*g$. We explain our notational conventions after this introduction. This document contains formulas such as the following:
\begin{thm*}
Let $X$ be a smooth section of $F^*TM$, for example the restriction of a vector field in $M$ to the image of $F$. The mean curvature vector and $\sigma(t)$-divergence of the tangential projection $X^T$ evolve by
\begin{align}
	-\partial_{t^i}\langle H, X^\perp\rangle &= -\langle L(\partial_{t^i}F^\perp), X^\perp\rangle -\langle H, \nabla_{t^i}X^\perp\rangle -\tr_{\sigma(t)}Rm(\cdot, \partial_{t^i}F^T,\cdot, X^\perp)\\ 
	&\qquad +\langle II(\partial_{t^i}F^T, \cdot), \nabla^\perp X^\perp\rangle_{N\Sigma_t\otimes T^*\Sigma_t }  -\langle \nabla^T\partial_{t^i}F^T, \nabla^TX^\perp\rangle_{T\Sigma_t\otimes T^*\Sigma_t}\\ 
	&\qquad -\Div_{\sigma(t)}\langle II(\partial_{t^i}F^T, \cdot), X^\perp\rangle_{N\Sigma_t}\\
	&=-\langle L(X^\perp), \partial_{t^i}F^\perp\rangle -\langle H, \nabla_{t^i}X^\perp\rangle-\tr_{\sigma(t)}Rm(\cdot, \partial_{t^i}F^T,\cdot, X^\perp)\\ 
	&\qquad +\Div_{\sigma(t)}\langle \nabla^\perp X^\perp, \partial_{t^i}F^\perp\rangle_{N\Sigma_t}-\Div_{\sigma(t)}\langle \nabla^\perp \partial_{t^i}F^\perp, X^\perp\rangle_{N\Sigma_t}\\
	&\qquad +\langle II(\partial_{t^i}F^T, \cdot), \nabla^\perp X^\perp\rangle_{N\Sigma_t\otimes T^*\Sigma_t } -\langle \nabla^T\partial_{t^i}F^T, \nabla^TX^\perp\rangle_{T\Sigma_t\otimes T^*\Sigma_t} \\
	&\qquad -\Div_{\sigma(t)}\langle II(\partial_{t^i}F^T, \cdot), X^\perp\rangle_{N\Sigma_t}\\
	+\partial_{t^i}(\Div_{\sigma(t)} X^T)&=	-\langle H, \nabla_{t^i}X^T\rangle +\Div_{\sigma(t)} \left(\nabla_{t^i}^TX^T\right)\\
	&\qquad -\left(\langle \nabla^T\partial_{t^i}F^\perp,\nabla^TX^T\rangle_{T\Sigma_t\otimes T^*\Sigma_t} +\tr \nabla^T_{\nabla^T\partial_{t^i}F^T}X^T\right) \\ 
	&\qquad -\tr_{\sigma(t)} \left(Rm(\cdot\,, \partial_{t^i}F^\perp, \cdot \, , X^T)+Rm(\cdot\,, \partial_{t^i}F^T, \cdot\,, X^T) \right) \\ 
	&\qquad +\langle II(\cdot, X^T), \nabla^\perp \partial_{t^i}F^\perp\rangle_{N\Sigma_t\otimes T^*\Sigma_t} +\langle II(\cdot, X^T), II(\cdot, (\partial_{t^i}F)^T)\rangle_{N\Sigma_t\otimes T^*\Sigma_t} 
\end{align}
Let $\sqrt{\sigma(t)}\in \Gamma(F^*(\Lambda^mT^*M))\rvert_\Sigma$ denote the induced (top-degree) area form on $\Sigma$ associated with $\sigma(t)$. We find
\begin{align}
	\partial^3_{t^it^jt^k}\sqrt{\sigma{(t)}}&=\mathcal T[\sigma(t), g, \nabla^\perp X, \cS{X}, Rm_g, \partial_iF, \nabla_i\partial_jF, \nabla_i\nabla_j\partial_kF]\sqrt{\sigma(t)}\\ 
	&\qquad +\partial_{t^i}\left(\sigma^{\alpha\beta}\cS{\partial_{t^k}F}_{\alpha\beta}\right)\partial_{t^j}\sqrt{\sigma(t)}+\partial_{t^j}\left(\sigma^{\alpha\beta}\cS{\partial_{t^k}F}_{\alpha\beta}\right)\partial_{t^i}\sqrt{\sigma(t)}\\ 
	&\qquad +\left(\sigma^{\alpha\beta}\cS{\partial_{t^k}F}_{\alpha\beta}\right)\partial^2_{t^it^j}\sqrt{\sigma(t)}
\end{align}
where $\cS{X}:=(\nabla^TX)^\flat\in \Gamma(T^*\Sigma_t^{\otimes 2})$ packages the ($\sigma(t)$-lowered) shape operator and tangential covariant derivative together in a computationally expedient manner and $\mathcal T$ comprises contractions of its tensorial arguments for (iterated covariant) parameter derivatives $X=\partial_iF, \nabla_i\partial_jF\ldots$, given explicitly in the main text. 
\end{thm*}
More identities (in particular, lower order derivatives of $\sqrt{\sigma(t)}$) are contained in the main body of the text, calculated in a way that makes further derivative computations more bearable. We then provide an example of their use in the minimal area rigidity literature. Further use will feature in forthcoming work.

\subsection{Motivation}

The formulas we choose to present are meant for the following  geometric inverse problem. On a (pseudo--)Riemannian manifold-with-boundary $(M^n,g)$, consider a $k$-parameter family of loops $\gamma_t: S^1\to \partial M$, $t=(t^1,\ldots, t^k)$. Presume there exists a unique $k$-parameter family $F:I^k\times \Sigma^2\to (M,g)$ of properly embedded minimal surfaces $F_t(\Sigma)$ with $\partial F_t(\Sigma)=\gamma_t$ whose areas $A(t):=Area(F_t(\Sigma))$ are known. Does the boundary geometry on $\partial M$ combined with knowledge of $A(t)$ for sufficiently many families $(\gamma_t)\subset \partial M$ uniquely determine $g$ up to isometry in the ambient interior region $\bigcup_{\gamma_t} F_t(\Sigma)$?

For Riemannian $g$, this is an analog of the classic boundary rigidity problem for geodesic lengths \cite{pestov2005two,stefanov2005boundary}. The study of the minimal area inverse problem for ambient dimension $n=3, M\approx \mathbb{B}^3, \Sigma\approx \mathbb{D}$ was pioneered in \cite{alexakis2020determining}, giving an affirmative answer to uniqueness presuming knowledge of $A(t)$ for perturbed families of foliating minimal disks and curvature assumptions on $g$. Recent work of \cite{alvarez2026foliated} concludes uniqueness for more general area-simple metrics $g$ from $A(t)$ corresponding to $3$-parameter families of $\gamma_t$. For general orientable Riemannian surfaces-with-boundary $\Sigma$, work of \cite{carstea2024calder} used higher order linearizations to first uniquely determine the isometry class o,f the induced metric $\sigma(0)$ of the minimal surface $F_0(\Sigma)$ using $A(t)$ corresponding to boundary curves $\gamma_t:=\exp^{N\Sigma_0}\rvert_{\gamma_0}(t^kf_k\vec{n})$, where $\vec{n}$ is the unit normal of $\Sigma_0$. The followup paper \cite{carstea2026recovering} then determines $g$ on foliated regions. 

Approaches to the codimension one minimal area rigidity problem thus far lie under two camps. Spiritually, the work of  \cite{alexakis2020determining} and \cite{carstea2024calder,carstea2026recovering} are similar in that they first recover (some information on) the induced metric $\sigma(t)$ via a Calder\'on-type problem for the stability operator $L$ on each $\Sigma(t)$, before ``pasting'' together such leaves to recover the bulk metric (by way of ``tilts'' in \cite{alexakis2020determining} and ``shifts'' in \cite{carstea2026recovering}). Higher order derivatives of mean curvature are a key part of these approaches, and have mostly been computed in a \textit{ad hoc} manner with particular choices of ambient coordinates.

Meanwhile, \cite{alvarez2026foliated} uses an integral geometry approach via a Radon-type transform for minimal surfaces, similar in spirit to \cite{ambrozio2025riemannian} (for a version of nonrigidity) in the closed minimal hypersurface case. Key to their analysis of minmal areas $A(t)$ is an identification of the space of boundary circles $\mathcal C\ni \gamma_t$ with their Gauss lifts in the ambient unit tangent bundle, recalling that for codimension one surfaces, $\sqrt{\sigma(t)}\in F^*(\Lambda^2T^*M)$ can be identified with an oriented normal field. An integral approach was also used in \cite{busch2025generalized}, which applied double-fibration transforms in a presumed analytic background to obtain rigidity from areas of minimal hypersurfaces $\Sigma^{n-1}$. See also \cite{grebnev2026linearized} for related work in a linearized setting in a Euclidean background.

While from a geometric point of view, it may be unnatural to consider situations where one knows minimal areas without some knowledge of the metric \textit{a priori}, a main thrust for such minimal area inverse problems comes from physics, where $\sigma(t)$ are Riemannian metrics induced on spacelike minimal (``extremal'') surfaces $F_t(\Sigma^2)$ of codimension two embedded in an ambient Lorentzian, asymptotically AdS (aAdS) background $g$ . In the context of the AdS/CFT correspondence \cite{maldacena1999large}, it was proposed in \cite{ryu2006aspects} and generalized by \cite{Hubeny:2007xt} that the `quantum entanglement entropy' of a region enclosed by the loop $\gamma_t$ lying in the conformal infinity $\partial M$ is given up to a factor by the (renormalized) area of the presumptive minimal $F_t(\Sigma^{int})\subset M^{int}$ with  $\partial F_t(\Sigma)=\gamma_t$. In the physics literature, the inverse problem of determining $g$ from renormalized $A(t)$ has been dealt with in works such as \cite{bao2019towards} by arguments analogous to \cite{alexakis2020determining}, but as the authors themselves note, their arguments are not mathematical proofs. As discussed in \cite[Ch. 13]{rangamani2017holographic}, further proposals of linking bulk and boundary operators require analysis of various other submanifold extrinsic geometries. We hope that our calculations contained below, without restriction to vanishing mean curvature in all normal directions, may find use in such applications.

Aspects of AdS/CFT have only recently been mathematically treated in \cite{guillarmou2025probabilistic}. There has been comparably more work in the purely geometric side of things. Recall that an asymptotically hyperbolic (AH) Riemannian metric $g$ arises as a spacelike restriction of an ambient static aAdS metric. The inverse problem of renormalized distance rigidity in AH setting was treated in \cite{graham2019x}. Variations of renormalized area for minimal submanifolds in the AH setting were treated by \cite{alexakis2010renormalized} in the codimension one surface case, then more generally in \cite{marx2021variations}. The latter features calculations including the second derivative of mean curvature corresponding to families $F:I^k_t\times I_x\times \Sigma\to (M, g)$ of graphical perturbations at $t=0$, where $x$ is a special boundary defining function for $\partial M$. This led to \cite{marx2024inverse} showing that renormalized minimal areas determine the conformal infinity on $\partial M$ as well as the expansion of $g$ in $x$ near $\partial M$.

\subsection{Outline}
We first fix our geometric setting and explain our approach. This is followed by computations of the pointwise evolution of mean curvature (and divergence). This then yields derivatives of induced area forms and, upon integration, those of minimal areas. We conclude with a worked example in the codimension one setting. 

\subsection{Acknowledgements}
The author thanks John M. Lee for looking over an early form of these calculations, as well as Tracey Balehowsky and Tony Liimatainen for useful discussions. A portion of this work was undertaken while the author was a postdoc of Spyros Alexakis at the University of Toronto.

As a postdoc of Gustav Holzegel, the author is funded by the Deutsche Forschungsgemeinschaft under Germany's Excellence Strategy EXC 2044/2 –390685587, Mathematics Münster: Dynamics–Geometry–Structure, as well as the Alexander von Humboldt Foundation in the framework of the Alexander
von Humboldt Professorship endowed by the Federal Ministry of Education and Research.

\section{Preliminaries}
We consolidate our notation for standard submanifold geometry: see for example \cite[Ch. 1]{colding2011course}, \cite[Ch. 6]{docarmo1992riemannian}, \cite[Ch. 8]{lee2018introduction}. In particular, we define and motivate our use of $\cS{X}$, which coincides with the (lowered) shape operator when $X$ is normal.
\subsection{Geometric setup}
Let $F:I^k\times \Sigma^m\to (M^n,g)$ be a $k$-parameter smooth family of isometric immersions as in the introduction. One pulls back $TM$ with $F$ to form the bundle $F^*(TM)\to (I^k\times \Sigma)$, and differentiates sections along curves in $I^k\times \Sigma $ using the Levi-Civita connection $\nabla$ of $g$. As an example, for $X=X^K(t,p)\partial_K\in \Gamma(F^*(TM))$ for $(x^K)$ local coordinates on $M$, we have $\nabla_{t^i}X=(\partial_{t^i}X^K+(\partial_{t^i}F^I)X^J\Gamma_{IJ}^K)\partial_K$ for $\nabla$-Christoffel symbols $\Gamma_{IJ}^K$.

Over each constant parameter slice $\{t\,\}\times \Sigma$, $t=(t^1, \ldots, t^k)$, the bundle splits orthogonally with respect to $g$ into the direct sum $F^*(TM)=T\Sigma\oplus N\Sigma_{t}$ of tangent and normal bundles, equipped with metrics $\sigma(t):=F_t^*g$ and $\theta(t)$ induced by $g$, along with compatible connections $\nabla^T, \nabla^\perp$ obtained by taking the respective projections of $\nabla$. We denote $\Sigma$ equipped with the metric $\sigma(t)$ by $\Sigma_t$ (and $T\Sigma$ equipped with said metric by $T\Sigma_t$), and by a mild abuse of notation (and terminology) we denote its area (induced volume) form by $\sqrt{\sigma(t)}$. 

 The vector-valued second fundamental form $II_t\in \Gamma(\Sigma_t; T^*\Sigma^{\otimes 2}\otimes N\Sigma_t)$ is given by  $II_t(X,Y):=(\nabla_XY)^\perp$ for $X,Y\in T\Sigma $, and its $\sigma(t)$-trace yields the vector-valued mean curvature $H_t:=\tr_{\sigma(t)} II_t$. Acting on sections $N\in \Gamma (N\Sigma_t)$, the Simons' operator is $\tilde A(N):= \sum_{\alpha, \beta}\langle II_t(e_\alpha, e_\beta), N\rangle_g II_t(e_\alpha, e_\beta)$ for a $\sigma(t)$- orthonormal basis $(e_\alpha)$ at $T_p\Sigma_t$. The stability operator on $N\Sigma_t$ is given by $L(N):=\Delta^\perp N +\tr_{\sigma(t)}(R_g(\cdot, N)\, \cdot)+\tilde A(N)$, where $\Delta^\perp:=-\left(\nabla^\perp\right)^*\nabla^\perp =+\tr_{\sigma(t)}(\nabla^\perp)^2$ is the normal bundle Laplacian and our convention for the ambient Riemann curvature is $R_g(X,Y)W:=\nabla_Y\nabla_XW-\nabla_X\nabla_YW+\nabla_{[X,Y]}W$.

 \subsection{Defining $\cS{X}$} 

As will become apparent in subsequent calculations, for $X\in \Gamma(F^*TM)$, we find it useful to define $\cS{X}:=(\nabla^TX)^\flat\in \Gamma(T^*\Sigma_t^{\otimes 2})$. In local coordinates $(x^\alpha)$ on $\Sigma$, we have $\cS{X}_{\alpha\beta}=\langle \nabla_\alpha X, \partial_\beta F\rangle_g=g_{IJ}X^I_{;\alpha}\partial_\beta F^J$.

Our utility for this comes from the following observation. Note that by linearity, for $V_t, W_t\in T\Sigma_t$, we have 
\begin{align}
	\cS{X}(V_t, W_t)=\cS{X^\perp}(V_t, W_t)+\cS{X^T}(V_t, W_t)=-\langle II_t(V_t, W_t), X^\perp\rangle_{\theta(t)}+\langle \nabla^T_{V_t}X^T, W_t\rangle_{\sigma(t)}
\end{align}
and by taking the $\sigma(t)$-trace
\begin{align}
	\sigma^{\alpha\beta}\cS{X}_{\alpha\beta}=-\langle H, X\rangle_g+\Div_{\sigma(t)}X^T\,.
\end{align}
By the well-known fact (or see our calculations below) that 
\begin{align}\partial_{t^i}\sqrt{\sigma(t)}=\sigma^{\alpha\beta}\cS{\partial_{t^i}F}_{\alpha\beta}\sqrt{\sigma(t)}
\end{align} and parameter derivatives are absorbed as 
\begin{align}\partial_{t^i}\cS{X}_{\alpha\beta}=\cS{\nabla_{t^i}X}_{\alpha\beta}+ (\text{remainder}),
\end{align}
 higher derivatives of area come from differentiating under the integral and factoring as 
 \begin{align}
 	\partial^n A(t)=\int_{\Sigma}\partial^n\sqrt{\sigma(t)}=\int_{\Sigma}\partial^{n-1}(\sigma^{\alpha\beta}\cS{\partial_tF}_{\alpha\beta}\sqrt{\sigma})=\sum_j\int_\Sigma \binom{n-1}{j}\partial^{n-1-j} \left(\sigma^{\alpha\beta}\cS{\partial_tF}_{\alpha\beta}\right)(\partial^j\sqrt{\sigma})\,.
 \end{align}
Further distributing the right hand side parameter derivatives on $\sigma^{\alpha\beta}\cS{\partial_tF}$ and absorbing them within $\cS{\cdot}$ yields boundary integrals from the tangential divergence term; vanishing bulk integrals due to (presumed) minimality; and remainder bulk integrals featuring terms like the ambient curvature. This decomposition is of use in inverse problems: we provide an example at the end of this paper highlighting this in \cite{carstea2024calder}. 

\subsection{Notational (ab)uses}
The parameters are denoted by lower case Latin indices, while Greek indices are with respect to arbirary coordinates on $\Sigma$. Note carefully that $(x^\alpha)$ should not be construed as coordinates in the ambient manifold $M$ and will always refer to local coordinates on the fixed smooth manifold $\Sigma$ when used for calculations. When $T$ is a $(0,2)$-tensor, notation such as $\sigma^{\alpha\beta}T_{\alpha\beta}$ is the $\sigma(t)$-trace of a $(0,2)$-tensor $T$ on $\Sigma$ and is $\Sigma$-coordinate independent. The symmetrization of $T_{\alpha\beta}$ is $T_{(\alpha\beta)}:=\frac{1}{2}\left(T_{\alpha\beta}+T_{\beta\alpha}\right)$.

The unlabelled angle bracket $\langle \cdot, \cdot\rangle $ is with respect to the ambient metric $g$. Contextually, this may end up coinciding with either tangential or normal metrics, depending on the arguments. The $g$-lowered curvature is $Rm(X,Y,W, Z):=\langle R(X,Y)W, Z\rangle$.

We will sometimes use interchangeably $\partial_\alpha F= dF(\partial_\alpha)=\alpha$ and the same for parameter derivatives. In taking covariant derivatives, we sometimes use $\nabla_{(\cdot)}^\perp X$ to denote $(\nabla_{(\cdot)}X)^\perp$; this coincides with the usage of $\nabla^\perp$ as a connection on the normal bundle only when the appropriate arguments lie in the correct spaces. As shorthand, when taking derivatives in the parameters $I^k$, we often use $t, s\ldots$ as individual parameters in $I^k$ instead of $t^1, t^2\ldots $ and abbreviate $\partial_{t^i}$ as $\partial_t, \partial_s, \ldots $. When not appearing as a derivative, $t$ will still be retained as a label for the full paramater $k$-tuple. For example, $\Sigma_t$ is the manifold $\Sigma$ with induced metric $\sigma(t)$, while $\partial_s\partial_t\sqrt{\sigma(t)}$ measures the second derivative of the induced area form for $\sigma(t)$ in some parameter factors $t^i, t^j$ relabelled as $s, t$. The $\sigma(t)$ induced metric on $\partial\Sigma$ is denoted $\sigma^\partial (t)$.

\section{Pointwise evolution of mean curvature}
Recall first the following preliminary pointwise identities which feature in our subsequent calculations.

\begin{lem}\label{lem:preliminary} Let $t$ denote one parameter in the smooth family of isometric immersions above, and let $X, Y$ be sections of $F^*TM$. With respect to any local coordinates $(x^\alpha)$ on $\Sigma$, the following hold pointwise:
\begin{align}
	\partial_t\sigma_{\alpha\beta}&=\cS{\partial_tF}_{\alpha\beta}+\cS{\partial_tF}_{\beta\alpha}\\
	\partial_t\sigma^{\alpha\beta}&=-\sigma^{\alpha\mu}\sigma^{\beta\nu}(\cS{\partial_tF}_{\mu\nu}+\cS{\partial_tF}_{\nu\mu})\\
	\partial_t \left(Rm(\alpha, X, \beta, Y)\right)&=(\nabla_tRm)(\alpha, X, \beta, Y)\\ 
	&\qquad +Rm(\nabla_\alpha t, X, \beta, Y)+Rm(\alpha, \nabla_t X, \beta, Y) \\ 
	&\qquad +Rm(\alpha, X, \nabla_\beta t, Y)+Rm(\alpha, X,  \beta, \nabla_tY) \\ 
	\partial_t\langle \nabla_\alpha X, \nabla_\beta Y\rangle &=\langle \nabla_\alpha\nabla_tX, \nabla_\beta Y\rangle +\langle \nabla_\alpha X, \nabla_\beta\nabla_tY\rangle \\ 
	&\qquad +Rm(\alpha, t, X, \nabla_\beta Y)+Rm(\beta, t, Y, \nabla_\alpha X)\\ 
	\partial_{t}\cS{X}_{\alpha\beta}&=\cS{\nabla_tX}_{\alpha\beta}-Rm(\alpha, \partial_tF, \beta, X) +\langle \nabla_\alpha X, \nabla_\beta \partial_tF\rangle
\end{align}
\end{lem}
\begin{proof}
These are standard computations that we briefly illustrate for completeness and to gain familiarity with our notation. The first line follows from $\partial_t\sigma_{\alpha\beta}=\partial_t\langle \partial_\alpha F, \partial_\beta F\rangle =\langle \nabla_t\partial_\alpha F,\partial_\beta F\rangle+\langle \partial_\alpha F,\nabla_t\partial_\beta F\rangle$ and using $[\partial_t F, \partial_\alpha F]=dF[\partial_t, \partial_\alpha]=0$ to interchange covariant derivatives. The second exploits $\partial_t(\sigma^{\alpha\beta }\sigma_{\beta\gamma})=0$ to find $\partial_t\sigma^{\alpha \beta}=-\left(\partial_t\sigma_{\gamma\delta}\right)\sigma^{\alpha\gamma}\sigma^{\beta\delta}$. The remaining identities follow from standard properties of the Levi-Civita connection $\nabla$ on tensor products. 
\end{proof}

The following lemmas will be used in our formulation of general variations of mean curvature.

\begin{lem}\label{lem:prelim1}
	For $F$ as above and $X, Y$ smooth sections of $F^*TM$, we have at any $t$:
	\begin{align}
		\langle \nabla^\perp Y^\perp, \nabla^\perp X^\perp\rangle_{N\Sigma_t\otimes T\Sigma_t} = \Div_{\sigma(t)}(\langle \nabla^\perp Y^\perp, X^\perp\rangle) -\langle \Delta^\perp Y^\perp, X^\perp\rangle 
	\end{align}
\end{lem}
\begin{proof}
Fix $t$. Define $\omega(\cdot):=\langle \nabla_{(\cdot)} Y^\perp, X^\perp\rangle \in T^*\Sigma_t$. For $V_t, W_t\in T\Sigma_t$, we have 
\begin{align}
(\nabla^T\omega)(V_t, W_t)&=W_t\langle \nabla_{V_t} Y^\perp, X^\perp\rangle -\langle \nabla_{\nabla^T_{W_t}V_t}Y^\perp, X^\perp\rangle \\ 
&=\langle \nabla_{W_t}\nabla^\perp _{V_t} Y^\perp , X^\perp \rangle -\langle \nabla_{\nabla^T_{W_t}V_t}Y^\perp , X^\perp\rangle + \langle \nabla^\perp_{V_t}Y^\perp, \nabla_{W_t}X^\perp\rangle \,.
\end{align}
Choosing $V_t=W_t=e_\alpha$ for a local $\sigma(t)$-orthonormal frame $(e_\alpha)$ concludes the calculation.
\end{proof}

We also have the following formula relating divergence on the normal bundle to a tangential divergence term:
\begin{lem}\label{lem:divnormtotan}
	For $F$ as above, let $X$ be a section of $F^*TM$ and $Y$ be a section of $dF(T(I^k\times T\Sigma))\subset F^*TM$. We then have for all $t$
\begin{align}
	\Div_{\sigma(t)}\left(\langle \nabla^\perp Y, X^\perp\rangle \right)-\Div_{\sigma(t)} \langle [\,\cdot\,, Y], X^\perp\rangle =-\Div_{\sigma(t)} \left(\nabla_YX^\perp\right)^T
\end{align}

In particular, if $Y=\partial_tF$, we have
\begin{align}
	\Div_{\sigma(t)}\left(\langle \nabla^\perp \partial_tF, X^\perp\rangle \right) =-\Div_{\sigma(t)}\left(\nabla_tX^\perp\right)^T\,. 
\end{align}
\end{lem}
\begin{proof}
For any section $V_t\in T\Sigma_t\subset \Gamma(F^*TM)$, we have  $\langle V_t, X^\perp\rangle \equiv 0$, so $Y \langle V_t, X^\perp\rangle \equiv 0$ and therefore $\langle \nabla_{V_t}Y , X^\perp\rangle +\langle [Y, V_t], X^\perp\rangle =-\langle V_t, (\nabla_{Y }X^\perp)^T\rangle$. Hence, 

\begin{align}
		 -W_t\langle \nabla_{V_t}Y , X^\perp\rangle&=W_t\langle  (\nabla_{Y}X^\perp)^T, V_t\rangle+W_t\langle [V_t, Y], X^\perp\rangle\\ 
		 & =\langle  \nabla_{W_t}\nabla^T_{Y}X^\perp, V_t\rangle +\langle  \nabla_{Y}X^\perp, \nabla^T_{W_t}V_t\rangle+W_t\langle [V_t, Y], X^\perp\rangle\,.
\end{align}
We then have
\begin{align}
	-\langle  \nabla_{W_t}\nabla^T_{Y}X^\perp, V_t\rangle =W_t\langle \nabla_{V_t}Y , X^\perp\rangle+\langle  \nabla_{Y}X^\perp, \nabla^T_{W_t}V_t\rangle+W_t\langle [V_t, Y], X^\perp\rangle
\end{align}

For the right hand side, define $\omega(\cdot):= \langle \nabla^\perp_{\cdot} Y, X^\perp\rangle$ and $\eta(\cdot):=\langle [\, \cdot\, , Y], X^\perp\rangle$ observe

	\begin{align}
		(\nabla^T\omega )(V, W)&= W\left(\omega (V)\right)-\omega(\nabla^T_W V)\\
		&=W \langle \nabla^\perp_{V } Y, X^\perp\rangle-\langle \nabla^\perp_{\nabla^T_W V} Y, X^\perp\rangle\\ 
		&= W \langle \nabla^\perp_{V } Y, X^\perp\rangle-\langle \nabla_{Y} (\nabla^T_W V)_t, X^\perp\rangle-\langle [(\nabla^T_W V)_t, Y ], X^\perp \rangle \\ 
		&= W \langle \nabla^\perp_{V } Y, X^\perp\rangle+\langle \nabla_{Y}X^\perp,\nabla^T_W V\rangle -\langle [(\nabla^T_W V)_t, Y ], X^\perp \rangle
	\end{align}

We then have
\begin{align}
	-\left(\nabla^T(\nabla^T_YX^\perp) \right)^\flat(V_t, W_t)&=(\nabla^T\omega)(V_t,W_t)-W_t\langle [V_t, Y], X^\perp\rangle+\langle [(\nabla^T_W V)_t, Y ], X^\perp \rangle\\
	&= (\nabla^T\omega)(V_t,W_t)-(\nabla^T \eta)(V_t, W_t)
\end{align}

At any $t$, setting $V_t=W_t=e_\alpha$ for a local $\sigma(t)$-orthonormal frame $(e_\alpha)$ concludes the main calculation. For $Y=\partial_tF$, note that
$[\partial_\alpha F, \partial_tF]=dF[\partial_\alpha, \partial_t]=0$ for any $\sigma(t)$-normal coordinates $(x^\alpha)$. 

\end{proof}

These lead us to the vector mean curvature evolution equations.

\begin{prop}[First derivative of mean curvature and divergence]
	Let $F:I^K\times \Sigma \to (M,g)$ be a smooth $k$-parameter family of isometric immersions of $\Sigma$ and $X$ a section of $F^*TM$. Let $t$ be a  parameter in any one of the factors of $I^k$. We then have 
\begin{align}
	-\partial_t\langle H, X^\perp\rangle &= -\langle L(\partial_tF^\perp), X^\perp\rangle -\langle H, \nabla_tX^\perp\rangle -\tr_{\sigma(t)}Rm(\cdot, \partial_tF^T,\cdot, X^\perp)\\ 
	&\qquad +\langle II(\partial_tF^T, \cdot), \nabla^\perp X^\perp\rangle_{N\Sigma_t\otimes T^*\Sigma_t }  -\langle \nabla^T\partial_tF^T, \nabla^TX^\perp\rangle_{T\Sigma_t\otimes T^*\Sigma_t} \\ 
	&\qquad -\Div_{\sigma(t)}\langle II(\partial_tF^T, \cdot), X^\perp\rangle_{N\Sigma_t} \\
	&=-\langle L(X^\perp), \partial_tF^\perp\rangle -\langle H, \nabla_tX^\perp\rangle-\tr_{\sigma(t)}Rm(\cdot, \partial_tF^T,\cdot, X^\perp)\\ 
	&\qquad +\Div_{\sigma(t)}\langle \nabla^\perp X^\perp, \partial_tF^\perp\rangle_{N\Sigma_t}-\Div_{\sigma(t)}\langle \nabla^\perp \partial_tF^\perp, X^\perp\rangle_{N\Sigma_t}\\
	&\qquad +\langle II(\partial_tF^T, \cdot), \nabla^\perp X^\perp\rangle_{N\Sigma_t\otimes T^*\Sigma_t } -\langle \nabla^T\partial_tF^T, \nabla^TX^\perp\rangle_{T\Sigma_t\otimes T^*\Sigma_t} \\
	&\qquad -\Div_{\sigma(t)}\langle II(\partial_tF^T, \cdot), X^\perp\rangle_{N\Sigma_t} \\
	+\partial_t(\Div_\sigma X^T)&=	-\langle H, \nabla_tX^T\rangle +\Div_{\sigma(t)} \left(\nabla_t^TX^T\right)\\
	&\qquad -\left(\langle \nabla^T\partial_tF^\perp,\nabla^TX^T\rangle_{T\Sigma_t\otimes T^*\Sigma_t} +\tr \nabla^T_{\nabla^T\partial_tF^T}X^T\right) \\ 
	&\qquad -\tr_{\sigma(t)} \left(Rm(\cdot\,, \partial_tF^\perp, \cdot \, , X^T)+Rm(\cdot\,, \partial_tF^T, \cdot\,, X^T) \right) \\ 
	&\qquad +\langle II(\cdot, X^T), \nabla^\perp \partial_tF^\perp\rangle_{N\Sigma_t\otimes T^*\Sigma_t} +\langle II(\cdot, X^T), II(\cdot, (\partial_tF)^T)\rangle_{N\Sigma_t\otimes T^*\Sigma_t} 
\end{align}

For later purposes, it is useful to record this together as 
\begin{align}
	\partial_t(\sigma^{\alpha\beta}\cS{X}_{\alpha\beta})&=-2\sigma^{\alpha\mu}\sigma^{\beta\nu}\cS{\partial_tF}_{(\alpha\beta)}\cS{X}_{\mu\nu}\\ 
	&\qquad +\sigma^{\alpha\beta}\left(\cS{\nabla_tX}_{\alpha\beta}-Rm(\alpha, \partial_tF, \beta, X) +\langle \nabla_\alpha X, \nabla_\beta \partial_tF\rangle \right)
\end{align}
\end{prop}

\begin{proof}
We use the standard product rule and the preliminary identities to obtain the expression for $\partial_t(\sigma^{\alpha\beta}\cS{X}_{\alpha\beta})$, then further expand
\begin{align}
	\partial_t(\sigma^{\alpha\beta}\cS{X}_{\alpha\beta})&=-\sigma^{\alpha\mu}\sigma^{\beta\nu}\cS{\partial_tF}_{\beta\alpha}\cS{X}_{\mu\nu}\\ 
	&\qquad +\sigma^{\alpha\beta}\left(\cS{\nabla_tX}_{\alpha\beta}-Rm(\alpha, \partial_tF, \beta, X) +\langle \nabla^\perp_\alpha X, \nabla^\perp_\beta \partial_tF\rangle \right)\,. 
\end{align}

By linearity of $\cS{\cdot}$, we have

\begin{align}
	-\partial_t\langle H, X^\perp\rangle &= -\sigma^{\alpha\mu}\sigma^{\beta\nu}\cS{\partial_tF}_{\beta\alpha}\cS{X^\perp}_{\mu\nu}\\ 
	&\qquad +\sigma^{\alpha\beta}\left(\cS{\nabla_tX^\perp}_{\alpha\beta}-Rm(\alpha, t, \beta, X^\perp) +\langle \nabla^\perp_\alpha X^\perp, \nabla^\perp_\beta \partial_tF\rangle \right)\\ 
	+\partial_t(\Div_\sigma X^T)&= -\sigma^{\alpha\mu}\sigma^{\beta\nu}\cS{\partial_tF}_{\beta\alpha}\cS{X^T}_{\mu\nu}-\sigma^{\alpha\beta}Rm(\alpha, \partial_tF, \beta, X^T) \\ 
	&\qquad +\sigma^{\alpha\beta}\left(\cS{\nabla_tX^T}_{\alpha\beta}+\langle \nabla^\perp_\alpha X^T, \nabla^\perp_\beta \partial_tF\rangle \right)
\end{align}

Repackaging the latter in terms of tangential and normal projections yields the expression for the derivative of divergence. For the variation of mean curvature, we further apply Lemmas \ref{lem:prelim1} and \ref{lem:divnormtotan} to find either 
\begin{align}
\sigma^{\alpha\beta}\langle \nabla^\perp_{\alpha} X^\perp, \nabla^\perp_\beta\partial_tF\rangle &= -\langle \Delta^\perp (\partial_tF)^\perp, X^\perp\rangle - \Div_{\sigma(t)}(\nabla_tX^\perp)^T-\Div_{\sigma(t)}\langle II(\cdot, \partial_tF^T), X^\perp\rangle + \langle II(\partial_tF^T, \, \cdot\,), \nabla^\perp X^\perp\rangle \\ 
&=-\langle \Delta^\perp X^\perp, \partial_tF^\perp\rangle+ \Div_{\sigma(t)}\langle \nabla^\perp X^\perp, \partial_tF^\perp\rangle + \langle II(\partial_tF^T, \, \cdot\,), \nabla^\perp X^\perp\rangle \,,
\end{align}
and therefore either 
\begin{align}
-\partial_t\langle H, X^\perp\rangle &= -\langle \nabla^T\partial_tF^T, \nabla^TX^\perp\rangle-\langle H, \nabla_tX^\perp\rangle -\tr_{\sigma(t)}Rm(\cdot, \partial_tF^T,\cdot, X^\perp)\\ 
&\qquad -\Div_{\sigma(t)}\langle II(\cdot, \partial_tF^T), X^\perp\rangle +\langle II(\partial_tF^T, \cdot), \nabla^\perp X^\perp\rangle -\langle L(\partial_tF^\perp), X^\perp\rangle  \\ 
&=-\langle \nabla^T\partial_tF^T, \nabla^TX^\perp\rangle-\langle H, \nabla_tX^\perp\rangle-\Div_{\sigma(t)}\langle II(\partial_tF^T, \cdot), X^\perp\rangle  -\Div_{\sigma(t)}\langle \nabla^\perp \partial_tF^\perp, X^\perp\rangle\\ 
&\qquad -\tr_{\sigma(t)}Rm(\cdot, \partial_tF^T,\cdot, X^\perp)+\langle II(\partial_tF^T, \cdot), \nabla^\perp X^\perp\rangle -\langle L(X^\perp), \partial_tF^\perp\rangle +\Div_{\sigma(t)}\langle \nabla^\perp X^\perp, \partial_tF^\perp\rangle\,.
\end{align}
where in the latter equality we used Lemma \ref{lem:divnormtotan} to expand $\Div_{\sigma(t)}(\nabla_tX^\perp)^T$. 
\end{proof}

We proceed to calculate the second derivative of mean curvature and divergence. The following are straightforward secondary identities for bookkeeping purposes. 
\begin{lem}\label{lem:secondary} Let $t$ denote one parameter in the smooth family of isometric immersions above, and let $X, Y$ be sections of $F^*TM$. With respect to any local coordinates $(x^\alpha)$ on $\Sigma$, the following hold pointwise:
\begin{align}
	\partial_t(\sigma^{\alpha\beta}\cS{X}_{\alpha\beta})-\sigma^{\alpha\beta}\cS{\nabla_tX}_{\alpha\beta}&= \sigma^{\alpha\beta}\left(\langle \nabla_\alpha X, \nabla_\beta \partial_tF\rangle-Rm( \beta, X,\alpha, \partial_tF) \right)\\ 
	&\qquad -2\sigma^{\alpha\mu}\sigma^{\beta\nu}\cS{\partial_tF}_{(\alpha\beta)}\cS{X}_{\mu\nu}\\ 
	\partial_t\left(\langle \nabla_\alpha X, \nabla_\beta Y\rangle -Rm(\alpha, X, \beta, Y)\right)&=\left(\langle \nabla_\alpha\nabla_tX, \nabla_\beta Y\rangle -Rm(\alpha, \nabla_tX, \beta, Y)\right)\\ 
	&\qquad +\left(\langle \nabla_\beta\nabla_tY, \nabla_\alpha X \rangle - Rm(\beta, \nabla_tY, \alpha, X) \right)\\
	&\qquad -Rm(\nabla_\beta Y, X, \alpha, t)-Rm(\nabla_\alpha X, Y, \beta, t) \\
	&\qquad -Rm(\nabla_\alpha t, X, \beta, Y)-Rm(\nabla_\beta t, Y, \alpha, X)-(\nabla_tRm)(\alpha, X, \beta, Y)
\end{align}
\end{lem}

This aids in the calculation of 
\begin{prop}[Second derivative of mean curvature and divergence]
Let $F:I^k\times \Sigma \to (M,g)$ be a smooth $k$-parameter family of isometric immersions of $\Sigma$ and $X$ a section of $F^*TM$. Let $s,t$ be parameters in any of the factors of $I^k$. We then have:
\begin{align}
	\partial_{st}^2(\sigma^{\alpha\beta}\cS{X}_{\alpha\beta})	&=\partial_s\left(\sigma^{\alpha\beta}\cS{\nabla_tX}_{\alpha\beta}\right)\\
	&\qquad +\partial_t \left(\sigma^{\alpha\beta}\cS{\nabla_sX}_{\alpha\beta}\right)-\sigma^{\alpha\beta}\cS{\nabla_t\nabla_sX}_{\alpha\beta}\\
	&\qquad -2\sigma^{\alpha\mu}\sigma^{\beta\nu}\cS{X}_{(\alpha\beta)}\cS{\nabla_s\partial_tF}_{\mu\nu}+\sigma^{\alpha\beta}\left(\langle \nabla_\alpha\nabla_s t, \nabla_\beta X\rangle -Rm(\alpha, \nabla_s\partial_tF, \beta, X)\right)\\  
	&\qquad -2\sigma^{ \alpha\mu}\sigma^{\beta\nu}\left(\cS{\partial_sF}_{(\mu\nu)}\langle \nabla_\alpha \partial_tF, \nabla_\beta X\rangle+\cS{X}_{(\mu\nu)}\langle \nabla_\alpha \partial_tF, \nabla_\beta \partial_sF\rangle+\cS{\partial_tF}_{(\mu\nu)}\langle \nabla_\alpha X, \nabla_\beta \partial_sF\rangle\right) \\ 
	&\qquad +2\sigma^{\alpha\mu}\sigma^{\beta\nu}\left(\cS{\partial_sF}_{(\mu\nu)}Rm(\alpha, t,\beta, X)+\cS{X}_{(\mu\nu)}Rm(\alpha, \partial_sF, \beta, \partial_tF)+\cS{\partial_tF}_{(\mu\nu)}Rm(\alpha, \partial_sF, \beta, X) \right)\\ 
	&\qquad +\sigma^{\alpha\beta}\left(Rm(\alpha, s, t, \nabla_\beta X)+Rm(\beta, s, X,\nabla_\alpha t)\right)\\
	&\qquad -\sigma^{\alpha\beta}\left((\nabla_sRm)(\alpha, t,\beta, X) +Rm(\nabla_\alpha s, t,\beta, X)+Rm(\alpha, t,\nabla_\beta s, X)\right)\\ 
	&\qquad -2\partial_s\left(\sigma^{\alpha\mu}\sigma^{\beta\nu}\right)\cS{\partial_tF}_{(\mu\nu)}\cS{X}_{\alpha\beta} 
\end{align}
\end{prop}
\begin{proof}
We apply the product rule and implement Lemmas \ref{lem:preliminary} and \ref{lem:secondary} as follows:
\begin{align}
	\partial_{st}^2(\sigma^{\alpha\beta}\cS{X}_{\alpha\beta})&=\partial_s\left(\sigma^{\alpha\beta}\partial_t\cS{X}_{\alpha\beta}\right) +\partial_s\left((\partial_t\sigma^{\alpha\beta})\cS{X}_{\alpha\beta} \right)\\ 
	&=\partial_s\left(\sigma^{\alpha\beta}\cS{\nabla_tX}_{\alpha\beta}\right)+\partial_s\left(\sigma^{\alpha\beta}\left(\langle \nabla_\alpha t, \nabla_\beta X\rangle -Rm(\alpha, t,\beta, X)\right)\right) \\
	&\qquad -\partial_s\left(2\sigma^{\alpha\mu}\sigma^{\beta\nu}\cS{\partial_tF}_{(\mu\nu)}\cS{X}_{\alpha\beta} \right)\\ 
	&=\partial_s\left(\sigma^{\alpha\beta}\cS{\nabla_tX}_{\alpha\beta}\right)+\left(\partial_s\sigma^{\alpha\beta}\right)\left(\langle \nabla_\alpha t, \nabla_\beta X\rangle-Rm(\alpha, t,\beta, X)\right)\\
	&\qquad +\sigma^{\alpha\beta}\partial_s\left(\langle \nabla_\alpha t, \nabla_\beta X\rangle-Rm(\alpha, t,\beta, X) \right)\\  
	&\qquad -2\sigma^{\alpha\mu}\sigma^{\beta\nu}\cS{X}_{\alpha\beta}\partial_s\cS{\partial_tF}_{(\mu\nu)}-2\sigma^{\alpha\mu}\sigma^{\beta\nu}\cS{\partial_tF}_{(\mu\nu)}\partial_s\cS{X}_{\alpha\beta} \\
	&\qquad -2\partial_s\left(\sigma^{\alpha\mu}\sigma^{\beta\nu}\right)\cS{\partial_tF}_{(\mu\nu)}\cS{X}_{\alpha\beta} \\ 
	&=\partial_s\left(\sigma^{\alpha\beta}\cS{\nabla_tX}_{\alpha\beta}\right)+\left(\partial_s\sigma^{\alpha\beta}\right)\left(\langle \nabla_\alpha t, \nabla_\beta X\rangle-Rm(\alpha, t,\beta, X)\right)\\
	&\qquad +\sigma^{\alpha\beta}\left(\langle \nabla_\alpha\nabla_sX, \nabla_\beta t\rangle -Rm(\alpha, \nabla_sX, \beta, t)\right)\\  
	&\qquad  +\sigma^{\alpha\beta}\left(\langle \nabla_\beta\nabla_st, \nabla_\alpha X \rangle - Rm(\beta, \nabla_st, \alpha, X) \right)\\
	&\qquad -\sigma^{\alpha\beta}\left(Rm(\nabla_\beta t, X, \alpha, s)+Rm(\nabla_\alpha X, t, \beta, s) \right) \\
	&\qquad -\sigma^{\alpha\beta}\left(Rm(\nabla_\alpha s, X, \beta, t)+Rm(\nabla_\beta s, t, \alpha, X)+(\nabla_sRm)(\alpha, X, \beta, t)\right) \\
	&\qquad -2\sigma^{\alpha\mu}\sigma^{\beta\nu}\cS{X}_{\alpha\beta}\partial_s\cS{\partial_tF}_{(\mu\nu)}-2\sigma^{\alpha\mu}\sigma^{\beta\nu}\cS{\partial_tF}_{(\mu\nu)}\partial_s\cS{X}_{\alpha\beta}\\
	&\qquad -2\partial_s\left(\sigma^{\alpha\mu}\sigma^{\beta\nu}\right)\cS{\partial_tF}_{(\mu\nu)}\cS{X}_{\alpha\beta} \\ 
	&=\partial_s\left(\sigma^{\alpha\beta}\cS{\nabla_tX}_{\alpha\beta}\right)+\left(\partial_s\sigma^{\alpha\beta}\right)\left(\langle \nabla_\alpha t, \nabla_\beta X\rangle-Rm(\alpha, t,\beta, X)\right)\\
	&\qquad +\partial_t \left(\sigma^{\alpha\beta}\cS{\nabla_sX}_{\alpha\beta}\right)-\sigma^{\alpha\beta}\cS{\nabla_t\nabla_sX}_{\alpha\beta}\\  
	&\qquad  -2\sigma^{\alpha\mu}\sigma^{\beta\nu}\cS{X}_{(\alpha\beta)}\cS{\nabla_s\partial_tF}_{\mu\nu}+\sigma^{\alpha\beta}\left(\langle \nabla_\beta\nabla_st, \nabla_\alpha X \rangle - Rm(\beta, \nabla_st, \alpha, X) \right)\\
	&\qquad -\sigma^{\alpha\beta}\left(Rm(\nabla_\beta t, X, \alpha, s)+Rm(\nabla_\alpha X, t, \beta, s) \right)\\
	&\qquad -\sigma^{\alpha\beta}\left(Rm(\nabla_\alpha s, X, \beta, t)+Rm(\nabla_\beta s, t, \alpha, X)+(\nabla_sRm)(\alpha, X, \beta, t)\right)\\
	&\qquad +2\sigma^{\alpha\mu}\sigma^{\beta\nu}\cS{X}_{(\mu\nu)}(Rm(\alpha, \partial_sF, \beta, t) -\langle \nabla_\alpha t, \nabla_\beta \partial_sF\rangle) \\
	&\qquad +2\sigma^{\alpha\mu}\sigma^{\beta\nu}\cS{\partial_tF}_{(\mu\nu)}(Rm(\alpha, \partial_sF, \beta, X) -\langle \nabla_\alpha X, \nabla_\beta \partial_sF\rangle )\\
	&\qquad -2\partial_s\left(\sigma^{\alpha\mu}\sigma^{\beta\nu}\right)\cS{\partial_tF}_{(\mu\nu)}\cS{X}_{\alpha\beta}\, .
\end{align}
A regouping of terms leads to our result.
\end{proof}

For later purposes, the following tertiary identities are useful in combination with the above:
\begin{lem}\label{lem:tertiary}
	Let $s, t$ denote parameters in the smooth family of isometric immersions above, and let $X$ be a section of $F^*TM$. 

	If $X=\partial_rF$ where $r$ is another parameter of $F$, we have
\begin{align}
	-2\sigma^{\alpha\mu}\sigma^{\beta\nu}\cS{X}_{(\alpha\beta)}\cS{\nabla_s\partial_tF}_{\mu\nu}+\sigma^{\alpha\beta}\left(\langle \nabla_\alpha\nabla_s t, \nabla_\beta X\rangle -Rm(\alpha, \nabla_s\partial_tF, \beta, X)\right)=\\ 
	+\partial_r \left(\sigma^{\alpha\beta}\cS{\nabla_s\partial_tF}_{\alpha\beta}\right)-\sigma^{\alpha\beta}\cS{\nabla_r\nabla_s\partial_tF}_{\alpha\beta}\,.
\end{align}

Generally we can decompose $X,\nabla_st$ relative to $T\Sigma_t\oplus N\Sigma_t $ and write
\begin{align}
	-2\sigma^{\alpha\mu}\sigma^{\beta\nu}\cS{X}_{(\alpha\beta)}\cS{\nabla_s\partial_tF}_{\mu\nu}&+\sigma^{\alpha\beta}\left(\langle \nabla_\alpha\nabla_s t, \nabla_\beta X\rangle -Rm(\alpha, \nabla_s\partial_tF, \beta, X)\right)=\\ 
	\qquad &-\langle L(\nabla^\perp_st), X^\perp \rangle +\Div_\sigma \langle \nabla^\perp\nabla^\perp_st, X^\perp \rangle \\
	&+ \langle \nabla^\perp_\alpha(\nabla_s t), \nabla^\perp_\beta X^T\rangle+\langle \nabla^\perp_\alpha(\nabla_s t)^T, \nabla^\perp_\beta X^\perp\rangle \\
	&-\sigma^{\alpha\mu}\sigma^{\beta\nu}\cS{\nabla_st}_{\alpha\beta}\cS{X^T}_{\nu\mu}-\sigma^{\alpha\mu}\sigma^{\beta\nu}\cS{(\nabla_st)^T}_{\alpha\beta}\cS{X^\perp}_{\mu\nu}\\ 
	&-\sigma^{\alpha\beta}\left(Rm(\alpha, \nabla_s\partial_tF, \beta, X^T)+Rm(\alpha, (\nabla_s\partial_tF)^T, \beta, X^\perp)\right)
\end{align}
\end{lem}

\begin{proof}
	The former follows from the secondary identity, while the latter follows from expanding 
	\begin{align}
		-\sigma^{\alpha\mu}\sigma^{\beta\nu}\cS{X}_{\beta\alpha}\cS{\nabla_st}_{\mu\nu}+\sigma^{\alpha\beta}\left(\langle \nabla_\alpha^\perp\nabla_s t, \nabla^\perp_\beta X\rangle -Rm(\alpha, \nabla_s\partial_tF, \beta, X)\right)
	\end{align}
	and applying Lemma \ref{lem:prelim1} with $Y^\perp=(\nabla_st)^\perp$.
\end{proof}
\section{Derivatives of area forms and minimal areas}
\subsection{Derivatives of Induced Area Forms}
Upon substituting a parameter derivative $X=\partial_tF$ in the identities derived in the previous section, we can now more easily calculate the higher order derivatives of $\sqrt{\sigma(t)}$. We begin with the basic
\begin{prop}[First derivative of the area form]
Let $s,t$ denote parameters in the smooth family of isometric immersions above. We have the following ($\Sigma$-coordinate independent) identity

\begin{align}
	\partial_t\sqrt{\sigma(t)}&=\left(\sigma^{\alpha\beta}\cS{\partial_tF}_{\alpha\beta}\right)\sqrt{\sigma(t)}=\left(-\langle H, (\partial_tF)^\perp\rangle +\Div_{\sigma(t)}(\partial_tF)^T\right)\sqrt{\sigma(t)}
\end{align}
\end{prop}
\begin{proof}
This is classic but included for convenience. In any local coordinates on $\Sigma$, we have $\sqrt{\sigma(t)}=\sqrt{\det[\sigma(t)]_{\alpha\beta}}\,dx \left(=\|\partial_\alpha F_t\wedge \partial_\beta F_t\|_g\,dx\right)$. Near say $t=0$, recall
\begin{align}
\det[\sigma(t)]_{\alpha\beta}=\det\left([\sigma(0)]_{\alpha\beta}\right)\det \left(\delta_{\alpha\beta}+t[\sigma(0)]^{\alpha\beta}(\partial_t[\sigma (t)]_{\alpha\beta})\rvert_0+O(t^2)\right)
\end{align}
so we find
\begin{align}
\frac{d}{dt}\det\left([\sigma(t)]_{\alpha\beta}\right)=\sigma^{\alpha\beta}\partial_t\sigma_{\alpha\beta}\det\left([\sigma(t)]_{\alpha\beta}\right).
\end{align}
The earlier formula for derivative of the inverse metric and symmetry of $\sigma^{\alpha\beta}$ concludes the calculation.
\end{proof}

The following is then direct from the formula for first variation of mean curvature and divergence with $X=\partial_tF$:

\begin{prop}[Second derivative of the area form]
Let $s, t$ denote two (possibly non-distinct) parameters in the smooth family of isometric immersions above. We have the following ($\Sigma$-coordinate independent) identity
\begin{align}
	\partial_s\partial_t\sqrt{\sigma(t)}&=\partial_s\left(\sigma^{\alpha\beta}\cS{\partial_tF}_{\alpha\beta}\right)\sqrt{\sigma(t)}+\left(\sigma^{\alpha\beta}\cS{\partial_tF}_{\alpha\beta}\right)\partial_s\sqrt{\sigma(t)}\\ 
	&=\left(-\partial_s\langle H, \partial_tF\rangle+\partial_s(\Div_{\sigma(t)}\partial_tF^T)\right)\sqrt{\sigma(t)}+\left(\sigma^{\alpha\beta}\cS{\partial_tF}_{\alpha\beta}\right)\partial_s\sqrt{\sigma(t)}\,. 
\end{align}
\end{prop}
\begin{proof}
	This follows from applying the first derivative of the area form, the product rule, then finally the derivative of mean curvature and divergence with $X=\partial_tF$.
\end{proof}

\begin{prop}[Third derivative of the area form] Let $s, t, r$ denote three (possibly non-distinct) parameters in the smooth family of isometric immersions above. We have the following ($\Sigma$-coordinate independent) identity
\begin{align}
	\partial^3_{str}\sqrt{\sigma(t)}&=\Bigg[\partial_s\left(\sigma^{\alpha\beta}\cS{\nabla_t\partial_rF}_{\alpha\beta}\right)+\partial_t \left(\sigma^{\alpha\beta}\cS{\nabla_s\partial_rF}_{\alpha\beta}\right) +\partial_r\left(\sigma^{\alpha\beta}\cS{\nabla_s\partial_tF}_{\alpha\beta}\right)\\ 
	&\qquad -\sigma^{\alpha\beta}\cS{\nabla_t\nabla_s\partial_rF}_{\alpha\beta}-\sigma^{\alpha\beta}\cS{\nabla_r\nabla_s\partial_tF}_{\alpha\beta}\,\\
	&\qquad -2\sigma^{ \alpha\mu}\sigma^{\beta\nu}\left(\cS{\partial_sF}_{(\mu\nu)}\langle \nabla_\alpha \partial_tF, \nabla_\beta \partial_rF\rangle+\cS{\partial_rF}_{(\mu\nu)}\langle \nabla_\alpha \partial_tF, \nabla_\beta \partial_sF\rangle+\cS{\partial_tF}_{(\mu\nu)}\langle \nabla_\alpha \partial_rF, \nabla_\beta \partial_sF\rangle\right) \\ 
	&\qquad +2\sigma^{\alpha\mu}\sigma^{\beta\nu}\left(\cS{\partial_sF}_{(\mu\nu)}Rm(\alpha, t,\beta, \partial_rF)+\cS{\partial_rF}_{(\mu\nu)}Rm(\alpha, \partial_sF, \beta, \partial_tF)+\cS{\partial_tF}_{(\mu\nu)}Rm(\alpha, \partial_sF, \beta, \partial_rF) \right)\\ 
	&\qquad +\sigma^{\alpha\beta}\left(Rm(\alpha, s, t, \nabla_\beta r)+Rm(\beta, s, \partial_rF,\nabla_\alpha t)\right)\\
	&\qquad -\sigma^{\alpha\beta}\left((\nabla_sRm)(\alpha, t,\beta, \partial_rF) +Rm(\nabla_\alpha s, t,\beta, \partial_rF)+Rm(\alpha, t,\nabla_\beta s, \partial_rF)\right)\\ 
	&\qquad -2\partial_s\left(\sigma^{\alpha\mu}\sigma^{\beta\nu}\right)\cS{\partial_tF}_{(\mu\nu)}\cS{\partial_rF}_{\alpha\beta} \Bigg]\sqrt{\sigma(t)}\\ 
	&\qquad +\partial_t\left(\sigma^{\alpha\beta}\cS{\partial_rF}_{\alpha\beta}\right)\partial_s\sqrt{\sigma(t)}+\partial_s\left(\sigma^{\alpha\beta}\cS{\partial_rF}_{\alpha\beta}\right)\partial_t\sqrt{\sigma(t)}+\left(\sigma^{\alpha\beta}\cS{\partial_rF}_{\alpha\beta}\right)\partial^2_{st}\sqrt{\sigma(t)}
\end{align}
\end{prop}

\begin{proof}
Applying the product rule as before, we are left to calculate the term $\partial_s\partial_t\left(\sigma^{\alpha\beta}\cS{\partial_rF}_{\alpha\beta}\right)$. Apply the second derivative of mean curvature and divergence with $X=\partial_rF$ and apply the substitution in Lemma \ref{lem:tertiary}.
\end{proof}

\subsection{Derivatives of minimal areas}
The following are immediate from integrating the derivatives of the area form.

\begin{prop}[Second derivative of minimal areas]
Consider normal sections $V_k\in C^\infty(N\Sigma_0)$ for $k=1, 2$ with $L(V_1)=0$. Suppose $F:I^2\times \Sigma \to M$ is a smooth $2$-parameter family of isometric immersions $F_{t}(\Sigma)\subset M$ of $\Sigma$ with parameter $t:=(t^1,t^2 )$ such that 
\begin{enumerate}
	\item $(\Sigma, F_0^*g)=\Sigma_0$
	\item $F_{t}\rvert_{\partial\Sigma}=\exp^{N\Sigma_0}\rvert_{\partial\Sigma}(t^kV_k\rvert_{\partial\Sigma})$
	\item $\partial_{t^k}F\rvert_{0}=V_k$
	\item $H\left(F_{(t^1,0)}(\Sigma)\right)\equiv 0$.
\end{enumerate}
Then we have:
\begin{align}
	\partial_{t^1}\left(\int_{\partial\Sigma} \langle \partial_{t^2}F, \nu^{\partial\Sigma_t}\rangle \sqrt{\sigma^\partial (t)}\right)\Big\rvert_0=\partial^2_{t^1t^2}A\rvert_0=\int_{\partial\Sigma}\langle \nabla_\nu^\perp V_1, V_2\rangle \sqrt{\sigma^{\partial}(0)}
\end{align}
\end{prop}
\begin{proof}
The first equality is clear by direct integration. The second follows by observing $\partial_{t^1}\sqrt{\sigma(t)}\rvert_0=0$, 

\begin{align}
	\partial_{t^1}\left(\sigma^{\alpha\beta}\cS{\partial_{t^2}F^T}_{\alpha\beta}\right)=\Div_{\sigma(t)}\left(\nabla_{t^1}\partial_{t^2}F\right)^T-\Div_{\sigma(t)}\left(\nabla_{t^1}\partial_{t^2}F^\perp \right)^T\,, 
\end{align}
and applying Lemma \ref{lem:divnormtotan}.
\end{proof}

\begin{prop}[Third derivative of minimal areas]
	Consider normal sections $V_k\in \Gamma(N\Sigma_0)$ for $k=1, 2, 3$ which solve $L(V_k)=0$. Suppose $F:I^3\times \Sigma \to M$ is a smooth $3$-parameter family of isometric immersions $F_{t}(\Sigma)\subset M$, $t:=(t^1,t^2, t^3 )$, of $\Sigma$ such that 
\begin{enumerate}
	\item $(\Sigma, F_0^*g)=(\Sigma, \sigma)=\Sigma_0$
	\item $F_{t}\rvert_{\partial\Sigma}=\exp^{N\Sigma_0}\rvert_{\partial\Sigma}(t^kV_k\rvert_{\partial\Sigma})$
	\item $\partial_{t^k}F\rvert_{\Sigma_0}=V_k$
	\item $H\left(F_{(t^1,0,0)}(\Sigma)\right)\equiv 0$; $H\left(F_{(0,t^2,t^3)}(\Sigma)\right)\equiv 0$.
\end{enumerate}
Then we have:
\begin{align}
	\partial_{{t^3}{t^2}{t^1}}&A\rvert_0-\partial_{t^3}\left(\int_{\partial\Sigma}\langle \nabla_{t^2}(\partial_{t^1}F), \nu^{\partial\Sigma_t}\rangle \sqrt{\sigma^\partial (t)}\right)\Big\rvert_0\\ 
	&-\partial_{t^2}\left(\int_{\partial\Sigma}\langle \nabla_{t^3}(\partial_{t^1}F), \nu^{\partial\Sigma_t}\rangle \sqrt{\sigma^\partial (t)}\right)\Big\rvert_0-\partial_{t^1}\left(\int_{\partial\Sigma}\langle \nabla_{t^3}(\partial_{t^2}F), \nu^{\partial\Sigma_t}\rangle \sqrt{\sigma^\partial (t)}\right)\Big\rvert_0\\
	&\qquad +\int_{\partial\Sigma}\langle \nabla_{t^3}\nabla_{t^2}(\partial_{t^1}F), \nu^{\partial\Sigma_t}\rangle\sqrt{\sigma^\partial (0)}+\int_{\partial\Sigma}\langle \nabla_{t^2}\nabla_{t^3}(\partial_{t^1}F), \nu^{\partial\Sigma_t}\rangle\sqrt{\sigma^\partial (0)} \\ 
	&=\int_{\Sigma}  \left(\dddot{\mathcal{O}}(\nabla^\perp V, \nabla^\perp V, V)+\dddot{\mathcal{O}}(\nabla^\perp V, V, V)+\dddot{\mathcal{O}}(V, V, V)\right)\,\sqrt{\sigma (0)}
\end{align}
where 
\begin{align}
	\dddot{\mathcal{O}}&(\nabla^\perp V, \nabla^\perp V, V):=\\ 
	&-2\Big\langle \cS{V_3}, \langle \nabla^\perp V_1, \nabla^\perp V_2\rangle_{N\Sigma_0}\Big\rangle_{T^*\Sigma_0^{\otimes 2}}-2\Big\langle \cS{V_2}, \langle \nabla^\perp V_3, \nabla^\perp V_1\rangle_{N\Sigma_0}\Big\rangle_{T^*\Sigma_0^{\otimes 2}}-2\Big\langle \cS{V_1}, \langle \nabla^\perp V_2, \nabla^\perp V_3\rangle_{N\Sigma_0}\Big\rangle_{T^*\Sigma_0^{\otimes 2}}
\end{align}
\begin{align}
	\dddot{\mathcal{O}}(\nabla^\perp V, V, V)&:=\tr_{\sigma_0}\left(Rm(\cdot, V_3, V_2, \nabla_{\cdot}^\perp V_1)+Rm(\cdot, V_2, V_3, \nabla_{\cdot}^\perp V_1)+Rm(\cdot, V_3, V_1, \nabla_{\cdot}^\perp V_2)+Rm(\cdot, V_2, V_1, \nabla_{\cdot}^\perp V_3)\right)\\ 
	\dddot{\mathcal{O}}(V,V,V)&:=\Big\langle \cS{V_3},Rm(\cdot, V_2, \cdot, V_1) \Big\rangle_{T^*\Sigma_0^{\otimes 2}}+\Big\langle \cS{V_2},Rm(\cdot, V_3, \cdot, V_1) \Big\rangle_{T^*\Sigma_0^{\otimes 2}}\\ 
	&\qquad +2 \tr_{\sigma_0, 23} \tr_{\sigma_0, 45}\tr_{\sigma_0, 61}\left(\cS{V_1}\otimes  \cS{V_2}\otimes \cS{V_3} \right)-\left(\tr_{\sigma_0}(\nabla Rm)(\cdot, V_3, \cdot, V_2, V_1)\right)
\end{align}
\end{prop}
\begin{proof}
This follows directly from integrating the corresponding area form derivative formula, applying the divergence theorem for boundary terms, and applying the simplifications in the hypothesis, noting that $\cS{X^\perp}$ is symmetric.
\end{proof}

Further derivative identities may be computed using similar methods as needed.

\section{Example: Codimension one calculations}

As a computational example, we translate the higher order linearizations of \cite{carstea2024calder} using variations of the area form $\sqrt{\sigma(t)}$ directly, as well as translate the equations satisfied by the lapse and shift in \cite{carstea2026recovering} using our general formula for pointwise evolution of mean curvature. 

Consider Fermi coordinates $(x^\alpha, y)$ over a fixed minimal $F_0(\Sigma^2)$, and a $4$-parameter family of minimal variations $F:I^4\times \Sigma^2\to (M^3, g)$ given graphically by $F_t(x)=(x^\alpha, u_t(x))$ corresponding to minimal fill-ins of the boundary loops $\gamma_t=\exp^{N\Sigma_0}\rvert_{\gamma_0}(t^kf_k(x)\vec n(0))=(x^\alpha, t^kf_k(x^\alpha))$, where $f_k\in C^\infty(\partial\Sigma)$, $k=1,\ldots, 4$, and 
\begin{align}
	\vec n(t)=\frac{\partial_y-g^{\alpha\mu}(\partial_\mu u_t)\partial_\alpha }{\sqrt{1+\|du_t\|_{g(x, u_t)}^2}}
\end{align}
is the unit normal to each orientable $F_t(\Sigma)$. In our notation, we take $X=\vec n(t)$. We relate notation between the standard mean curvature equation, that of \cite{carstea2024calder}, and our formulation respectively by
\begin{align}
	\Div_g(\vec n)=-\frac{1}{\sqrt{ \det g_u}}\nabla\cdot \left(g_u^{-1}\frac{\sqrt{ \det g_u}}{\sqrt{ 1+|\nabla u|^2_{g_u}}}\right)\nabla u+f(u, \nabla u)=\sigma^{\alpha\beta}\cS{\vec n }_{\alpha\beta}
\end{align}

We stress that the right hand side is $\tr_{\sigma(t)}\cS{\vec{n}(t)}$, as the Greek indices appearing in the contraction are with respect to the frames $\partial_\alpha F=\partial_{x^\alpha}+\partial_\alpha u\partial_y$, but we keep this notation to be consistent with the previous section. For brevity in the following, notate $\partial_i:=\partial_{t^i},\, v_i(x):=\partial_{i}u\rvert_{t=0},\, v_{ij}(x):=\partial^2_{ij}u\rvert_{t=0}$. In these coordinates, $\partial_iF=(\partial_iu)\partial_y$. 

The first, second, and third linearizations of mean curvature in \cite{carstea2024calder} are obtained from calculating 
\begin{align}
	\partial_1\left(\sigma^{\alpha\beta}\cS{\vec{n}}_{\alpha\beta}\right)\Big\rvert_0,\, \partial^2_{21}\left(\sigma^{\alpha\beta}\cS{\vec{n}}_{\alpha\beta}\right)\Big\rvert_0, \,\partial^3_{321}\left(\sigma^{\alpha\beta}\cS{\vec{n}}_{\alpha\beta}\right)\Big\rvert_0
\end{align}
with respective integral identities given by integrating the $2$-forms 
\begin{align}
	v_1\,\partial_2\left(\sigma^{\alpha\beta}\cS{\vec{n}}_{\alpha\beta}\right)\Big\rvert_0\,\sqrt{\sigma(0)},\, v_1\,\partial^2_{32}\left(\sigma^{\alpha\beta}\cS{\vec{n}}_{\alpha\beta}\right)\Big\rvert_0\,\sqrt{\sigma(0)}, \, v_1\,\partial^3_{432}\left(\sigma^{\alpha\beta}\cS{\vec{n}}_{\alpha\beta}\right)\Big\rvert_0\,\sqrt{\sigma(0)}
\end{align}

This is obtained directly from area form variations as follows. Decompose $\partial_iF=V_i(t,x)\vec{n}(t)+(\partial_iF)^T$ and observe $V_i(t,x)\vec{n}\in N\Sigma_t$, hence $\cS{V_i\vec{n}}_{\alpha\beta}=V_i\cS{\vec{n}}$ at all parameters $t$. Under general possibly non-minimal variations, one finds 
\begin{align}
v_1\partial_2\left(\sigma^{\alpha\beta}\cS{\vec{n}}_{\alpha\beta}\right)\Big\rvert_0\,\sqrt{\sigma(0)}&=\partial^2_{21}\sqrt{\sigma}\Big\rvert_0-\partial_2\left(\sigma^{\alpha\beta}\cS{\partial_1F^T}_{\alpha\beta}\sqrt{\sigma}\right)\Big\rvert_{0}-(\partial_2V_1)\rvert_0\left(\sigma^{\alpha\beta}\cS{\vec{n}}_{\alpha\beta}\sqrt{\sigma}\right)\Big\rvert_{0}\\ 
v_1\,\partial^2_{32}\left(\sigma^{\alpha\beta}\cS{\vec{n}}_{\alpha\beta}\right)\Big\rvert_0\,\sqrt{\sigma(0)}&=\partial^3_{321}\sqrt{\sigma}\Big\rvert_0-\partial_{32}\left(\sigma^{\alpha\beta}\cS{\partial_1F^T}_{\alpha\beta}\sqrt{\sigma}\right)\Big\rvert_0 \\
&-\partial_3\left((\partial_2V_1)\left(\sigma^{\alpha\beta}\cS{\vec{n}}_{\alpha\beta}\sqrt{\sigma}\right)\right)\Big\rvert_{0}\\
&-(\partial_3V_1)\rvert_0\partial_2\left(\sigma^{\alpha\beta}\cS{\vec{n}}_{\alpha\beta}\right)\Big\rvert_{0}\sqrt{\sigma(0)}-v_1\partial_2\left(\sigma^{\alpha\beta}\cS{\vec{n}}_{\alpha\beta}\right)\Big\rvert_{0}\partial_3(\sqrt{\sigma})\rvert_0\\ 
v_1\,\partial^3_{432}\left(\sigma^{\alpha\beta}\cS{\vec{n}}_{\alpha\beta}\right)\Big\rvert_0\,\sqrt{\sigma(0)}&=\partial^4_{4321}\sqrt{\sigma}\Big\rvert_0-\partial^3_{432}\left(\sigma^{\alpha\beta}\cS{\partial_1F^T}_{\alpha\beta}\sqrt{\sigma}\right)\Big\rvert_0  \\ 
&\qquad -\partial^2_{43}\left((\partial_2V_1)\left(\sigma^{\alpha\beta}\cS{\vec{n}}_{\alpha\beta}\sqrt{\sigma}\right)\right)\Big\rvert_{0}\\
&\qquad -\partial_4\left((\partial_3V_1)\partial_2\left(\sigma^{\alpha\beta}\cS{\vec{n}}_{\alpha\beta}\right)\sqrt{\sigma}\right)\Big\rvert_{0}-\partial_4\left(V_1\partial_2\left(\sigma^{\alpha\beta}\cS{\vec{n}}_{\alpha\beta}\right)\partial_3(\sqrt{\sigma})\right)\Big\rvert_{0}\\
&\qquad -(\partial_4V_1)\rvert_0\partial^2_{32}\left(\sigma^{\alpha\beta}\cS{\vec{n}}_{\alpha\beta}\right)\Big\rvert_0\sqrt{\sigma(0)}-v_1\,\partial^2_{32}\left(\sigma^{\alpha\beta}\cS{\vec{n}}_{\alpha\beta}\right)\Big\rvert_0(\partial_4\sqrt{\sigma})\rvert_0
\end{align}
As in \cite{carstea2024calder} for the left hand side, we now add the assumptions that we are varying through minimals via initially normal variations, which results in $\partial^k\left(\sigma^{\alpha\beta}\cS{\vec{n}}_{\alpha\beta}\sqrt{\sigma}\right)\equiv 0$ and $\partial_i\sqrt{\sigma}\rvert_0=0$. These simplify the right hand sides to reveal the identities above as
\begin{align}
	v_1\partial_2\left(\sigma^{\alpha\beta}\cS{\vec{n}}_{\alpha\beta}\right)\Big\rvert_0\,\sqrt{\sigma(0)}&=\partial^2_{21}\sqrt{\sigma}\Big\rvert_0-\partial_2\left(\sigma^{\alpha\beta}\cS{\partial_1F^T}_{\alpha\beta}\right)\Big\rvert_{0}\sqrt{\sigma(0)}\\
	v_1\,\partial^2_{32}\left(\sigma^{\alpha\beta}\cS{\vec{n}}_{\alpha\beta}\right)\Big\rvert_0\,\sqrt{\sigma(0)}&=\partial^3_{321}\sqrt{\sigma}\Big\rvert_0-\partial^2_{32}\left(\sigma^{\alpha\beta}\cS{\partial_1F^T}_{\alpha\beta}\sqrt{\sigma}\right)\Big\rvert_0\\
	v_1\,\partial^2_{432}\left(\sigma^{\alpha\beta}\cS{\vec{n}}_{\alpha\beta}\right)\Big\rvert_0\,\sqrt{\sigma(0)}&=\partial^4_{4321}\sqrt{\sigma}\Big\rvert_0-\partial^3_{432}\left(\sigma^{\alpha\beta}\cS{\partial_1F^T}_{\alpha\beta}\sqrt{\sigma}\right)\Big\rvert_0
\end{align}
The integrals of the left hand side are the integral identities of \cite{carstea2024calder}. The right hand sides above integrate to become
\begin{align}
	0&=\partial^2_{21}A(t)\rvert_0-\int_{\partial\Sigma}v_1(\partial_\nu v_2)\, \sqrt{\sigma^\partial (0)}\\ 
	0&=\partial^3_{321}A(t)\rvert_0-\partial^2_{32}\left(\int_{\partial\Sigma}\langle \partial_1F, \nu^{\partial\Sigma(t)}\rangle_{\sigma(t)}\sqrt{\sigma^\partial  (t)}\right)\Big\rvert_0\\
	0&=\partial^4_{4321}A(t)\rvert_0-\partial^3_{432}\left(\int_{\partial\Sigma}\langle \partial_1F, \nu^{\partial\Sigma(t)}\rangle_{\sigma(t)}\sqrt{\sigma^\partial  (t)}\right)\Big\rvert_0
\end{align}

In codimension one, recall the standard identification of the normal bundle stability operator with a scalar elliptic operator by $L(\varphi):=\langle L(\varphi \vec n), \vec n\rangle $. The equation satisfied by the lapse derived in \cite{carstea2026recovering} for bulk metric $g$ determination is our equation for $\partial_t\langle H, X^\perp\rangle$ with $X^\perp=\vec{n}(t)$ and lapse $u:=\langle \partial_tF, \vec{n}\rangle_g $, while an equation for the ``shift'' vector field $(\partial_tF)^T$ can be derived from our equation for $\partial_t(\Div_{\sigma(t)}X^T)$  for $X^T:=(\partial_tF)^T$ upon fixing a choice of background coordinates, e.g. $H$-harmonic coordinates in \cite{carstea2026recovering}.

\bibliographystyle{plain} 
\bibliography{refs}

\end{document}